\documentclass[11pt,twoside,leqno]{aomamlt2e} 
  \pageno{725}
\received{June 20, 2001}

\newtheorem{lem}{Lemma}[section]
\newtheorem{prop}{Proposition}

\newtheorem{coro}{Corollary}

\newtheorem{fact}{Fact}
\newtheorem{theo}{Theorem}

\newenvironment{nproof}[1]{\trivlist\item[\hskip \labelsep{\bf Proof{#1}.}]}
{\hfill$\square$\endtrivlist}

\newenvironment{block}[1]{\trivlist\item[\hskip \labelsep{{#1}.}]}{\endtrivlist}

\theorembodyfont{\upshape}

\newtheorem{defi}{{\it Definition}}
 
\begin{document}
\currannalsline{159}{2004}

\def\itemm#1{ \item{\makebox[0.3in][l]{#1}}}

\newcommand{\Cr}{{\bf Cr}}
\newcommand{\dist}{{\rm dist}}
\newcommand{\diam}{\mbox{diam}\, }
\renewcommand{\mod}{{\rm mod}\,}
\newcommand{\compose}{\circ}
\newcommand{\dbar}{\bar{\partial}}
\newcommand{\Def}[1]{\emph{#1}}
\newcommand{\dx}[1]{\frac{\partial #1}{\partial x}}
\newcommand{\dy}[1]{\frac{\partial #1}{\partial y}}
\newcommand{\Res}[2]{{#1}\raisebox{-.4ex}{$\left|\,_{#2}\right.$}}

\newcommand{\C}{{\bf C}}
\newcommand{\D}{{\bf D}}
\newcommand{\Dm}{{\bf D_-}}
\newcommand{\N}{{\bf N}}
\newcommand{\R}{{\bf R}}
\newcommand{\Z}{{\bf Z}}

\font\mathfonta=msam10 at 11pt
\font\mathfontb=msbm10 at 11pt
\def\Bbb#1{\mbox{\mathfontb #1}}
\def\lesssim{\;\mbox{\mathfonta.}\;}
\def\suppset{\mbox{\mathfonta{c}}}
\def\subbset{\mbox{\mathfonta{b}}}
\def\grtsim{\;\mbox{\mathfonta\&}\;}
\def\gtrsim{\mbox{\mathfonta\&}}
\def\area{\mathrm{area}}
\def\itemm#1{ \item{\makebox[0.3in][l]{#1}}}
\def\diam{{\mathrm{diam}}}
\def\fr{\frac}
\def\eps{\varepsilon}
\def\om{\Omega}
\def\del{\Delta}
\def\al{\alpha}
\def\scal#1#2{\left\langle{#1},{#2}\right\rangle}
\def\br#1{\left(#1\right)}
\def\brr#1{\left[#1\right]}
\def\brs#1{\left\{#1\right\}}
\def\abs#1{\left|#1\right|}
\def\norm#1{\left\|#1\right\|}
\def\dd{\partial}
\def\Ai{A_\infty}
\def\J{J_F}
\def\spt{{\mathrm{supp}}}
\def\Var#1{{{\mathrm Var}\br{#1}}}
\def\qhdist#1{{{\mathrm{dist_{qh}}}\br{#1}}}
\def\qhl#1{{{\mathrm{l_{qhyp}}}\br{#1}}}
\def\const{{\mathrm{const}}}
\def\absconst{{\mathrm{abs.const.}}}
\def\Crit{{\mathrm{Crit}}}
\def\dg{{\mathrm{deg}}}
\def\HD{{\mathrm{HD}}}
\def\kerr{{\mathrm{ker}}}
\def\eon{{e^{o(n)}}}
\def\eonasymp{\overset{\eon}\to\asymp}
\def\m#1{\mu(#1)}
\def\mmax{\mu_{max}}
\def\supF{{\sup|F'|}}
\def\Cal#1{{\cal#1}}
\def\xx#1#2{{{\cal X}_{#1}(#2)}}
\def\diam{{\mathrm{diam}}}
\def\degg#1#2{{\# \br{#1,#2}}}

\def\len#1{{{\mathrm{length}}\br{#1}}}
\def\qhlen#1{{{\mathrm{length_{qh}}}\br{#1}}}

\def\typei{{\bf I}}
\def\typeip{{\bf I'}}
\def\typeii{{\bf II}}
\def\typeiii{{\bf III}}
\def\Sig{{\bf\Sigma}}

\font\mathfonta=msam10 at 11pt
\font\mathfontb=msbm10 at 11pt
\def\Bbb#1{\mbox{\mathfontb #1}}
\def\lesssim{\mbox{\mathfonta.}}
\def\suppset{\mbox{\mathfonta{c}}}
\def\subbset{\mbox{\mathfonta{b}}}
\def\grtsim{\mbox{\mathfonta\&}}
\def\gtrsim{\mbox{\mathfonta\&}}

\newcommand{\Bo}{\Box^{n}_{i}}
\newcommand{\co}{{\mbox con}}
\newcommand{\Di}{{\cal D}}
\newcommand{\gd}{{\underline \gamma}}
\newcommand{\gu}{{\underline g }}
\newcommand{\ce}{\mbox{III}}
\newcommand{\be}{\mbox{II}}
\newcommand{\F}{\cal{F_{\eta}}}
\newcommand{\Ci}{\bf{C}}
\newcommand{\ai}{\mbox{I}}
\newcommand{\dupap}{\partial^{+}}
\newcommand{\dm}{\partial^{-}}
\newenvironment{note}{\begin{sc}{\bf Note}}{\end{sc}}
\newenvironment{notes}{\begin{sc}{\bf Notes}\ \par\begin{enumerate}}%
{\end{enumerate}\end{sc}}
\newenvironment{sol}
{{\bf Solution:}\newline}{\begin{flushright}
{\bf QED}\end{flushright}}

\title{Metric attractors\\ for smooth unimodal maps}

 \acknowledgements{The third author was partially supported by NSF grant DMS-0072312.}
 
\twoauthors{Jacek Graczyk, Duncan Sands,}{Grzegorz \'{S}wia\c\negthinspace tek}

 \institution{University of Paris XI, Orsay, France\\
\email{Jacek.Graczyk@math.u-psud.fr}\\
\vglue-4pt
University of Paris XI, Orsay, France\\
\email{Duncan.Sands@math.u-psud.fr}\\ \vglue-4pt
Pennsylvania State University, University Park, State College, PA\\
\email{swiatek@math.psu.edu}
}

 \shorttitle{Metric attractors for smooth unimodal maps}
 \shortname{Jacek Graczyk, Duncan Sands, and Grzegorz \'{S}wia\c\negthinspace tek}

\centerline{\bf Abstract}
\vglue12pt

We classify the measure theoretic attractors of general $C^3$ unimodal maps
with quadratic critical points. The main ingredient is the decay of geometry.
 
\section{Introduction}

 1.1. {\it  Statement of results}.
The study of measure theoretical attractors occupied a central
position  in the theory of smooth dynamical systems in the 1990s.
Recall that a 
forward invariant compact set $A$ is called a (minimal)
{\em metric attractor} for some
dynamics if the basin of attraction  $B(A):=\break \{x : \omega(x)\subset A\}$  of
$A$ has positive Lebesgue measure and
$B(A')$ has Lebesgue measure zero for every forward
invariant compact set $A'$ strictly contained in $A$.
Recall that a set is \emph{nowhere dense} if its closure
has empty interior, and \emph{meager} if it is
a countable union of nowhere dense sets.
A forward invariant compact set $A$ is called a
(minimal) \Def{topological attractor} if $B(A)$ is not meager
while $B(A')$ is meager for every forward invariant
compact set $A'$ strictly contained in $A$.
A basic question, known as Milnor's problem, is whether
the metric and topological  attractors coincide for a 
given smooth unimodal map.  

Milnor's problem has a long and turbulent history;
see ~\cite{mest}, [11], [5], \cite{nesz}.
In the class of $C^{3}$ 
unimodal maps with negative Schwarzian derivative and a quadratic
critical point, an early solution to Milnor's problem
was given in~\cite{swia}. Recently, it was discovered that 
~\cite{swia} does 
not provide a complete proof. The author has told us that
his argument can
be repaired,~\cite{swiae}.
A correct solution using different techniques can be found in~\cite{nesz}.
A negative solution when the critical point  has high order 
is given in~\cite{bkns}.
The $C^{3}$ stability theorem of~\cite{kpd}, \cite{kozl1} implies that 
a generic $C^{3}$ unimodal map has finitely many metric 
attractors which are all attracting cycles, thus solving Milnor's problem in
the generic case.

Our current work solves Milnor's problem for smooth unimodal maps with
a quadratic critical point.
Historically, the solution is  based on two key developments.
The first, \cite{nesz}, established decay of geometry 
for a  class of $C^{3}$ nonrenormalizable
box mappings with finitely many branches
and negative Schwarzian derivative everywhere
except at the critical point which must be quadratic.
The second, \cite{kozl}, recovers negative Schwarzian derivative
for smooth unimodal maps with nonflat critical point:
the first return map to a neighborhood of the critical value has
negative Schwarzian derivative. 

Technically, for our study of metric attractors in the smooth category we need
a different estimate from that of~\cite{kozl}, one which works near
the critical point rather than the critical value~\cite{comp}.
We add a new Koebe lemma and exploit the fact that negative
Schwarzian derivative is not an invariant of smooth conjugacy to
show that the first return map to a neighborhood of the critical point
can be real-analytically conjugated to one having negative Schwarzian
derivative.
This makes it easy to transfer results known for maps with negative
Schwarzian to the smooth class.
Earlier results in this direction, in particular that high iterates
of a smooth critical circle homeomorphism have negative Schwarzian,
were obtained in~\cite{victim}. 

The classification of metric attractors containing the (nondegenerate)
critical point was announced in~\cite{comp}.  Here we give full proofs and 
explain the structure of metric attractors not containing the critical point
(based on the work of Ma\~n\'e~\cite{mane}). 
Consequently, we obtain the classification of all metric
attractors for smooth unimodal maps with a nondegenerate critical point.

\demo{Classification of metric dynamics}
A $C^1$ map $f$ of a compact interval $I$ is \Def{unimodal} if it has
exactly one point $\zeta$ where $f'(\zeta)=0$ (the \Def{critical point}),
$\zeta \in {\rm int}\,I$, $f'$ changes sign at $\zeta$,
and $f$ maps the boundary of $I$ into itself.
The critical point of $f$ is \Def{$C^n$ nonflat of order $\ell$} if,
near $\zeta$, $f$ can be written as
$f(x) = \pm |\phi(x)|^\ell + f(\zeta)$ where
$\phi$ is a $C^{n}$ diffeomorphism.
The critical point is \Def{$C^n$ nonflat} if it is $C^n$ nonflat of
some order $\ell > 1$. The set of  critical points of $f$ 
is denoted by $\Crit$.

\begin{theo}\label{tw2}
Let $I$ be a compact interval and $f : I \to I$
be a $C^{3}$ unimodal map with $C^{3}$ nonflat critical point of order $2$.
Then the $\omega$\/{\rm -}\/limit set of Lebesgue almost every point of $I$ is either
\begin{itemize}  \itemsep=2pt
\ritem{1.} a nonrepelling periodic orbit{\rm ,} or
\ritem{ 2.} a transitive cycle of intervals{\rm ,} or
\ritem{3.} a Cantor set of solenoid type. \end{itemize}
\end{theo}

A compact interval $J$ is \Def{restrictive} if $J$ contains the critical
point of $f$ in its interior and, for some $n > 0$, $f^n(J) \subseteq J$ and
$\Res{f^n}{J}$ is unimodal.
In particular, $f^n$ maps the boundary of
$J$ into itself.  This restriction of $f^n$ to $J$ is called a
\Def{renormalization} of $f$.
We say that $f$ is \Def{infinitely renormalizable}
if it has infinitely many restrictive intervals.

A periodic point $x$ of period $p$ is \Def{repelling} if
$|Df^p(x)| > 1$, \Def{attracting} if $|Df^p(x)| < 1$,
\Def{neutral} if $|Df^p(x)| = 1$ and \Def{super-attracting} if
$Df^p(x)=0$.  It is \Def{topologically attracting} if its
\Def{basin of attraction} $B(x) :=
B(\{x, f(x), \ldots, f^{p-1}(x)\})$ has nonempty interior.

A \Def{transitive cycle of intervals} is a finite union $C$ of
compact intervals such that $C$ is invariant under $f$, $C$ contains the
critical point of $f$ in its interior, and the action of $f$ on $C$ is
transitive (has a dense orbit).

We say that $f$ has a \Def{Cantor set of solenoid type}
if $f$ is infinitely renormalizable, the solenoid then being the
$\omega$-limit set of the critical point.

Note that the critical point $\zeta$ of a $C^4$ unimodal map with
$f''(\zeta)\neq 0$ is $C^3$ nonflat of order $2$.
The fact that the critical point has order 2 is used in an essential
way to exclude the possibility of \Def{wild Cantor attractors}.

\begin{coro}\label{coro:1}
Every metric attractor of $f$ is either 
\vglue-28pt
\phantom{ha}
\begin{itemize}\itemsep=2pt
\ritem{1.} a topologically attracting periodic orbit{\rm ,} or
\ritem{2.} a transitive cycle of intervals{\rm ,} or
\ritem{3.} a Cantor set of solenoid type.
\end{itemize}
\vglue-12pt\noindent
There is at most one metric attractor of type other than $1$.
\end{coro}
$$
 \epsfig{file=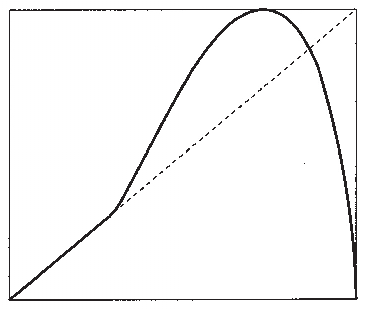}
$$
\centerline{Figure 1: Almost every point is mapped into the interval of fixed points.}

\vglue4pt
Figure~1 shows a unimodal map satisfying our hypotheses for
which the $\omega$-limit set of Lebesgue almost-every point is a
neutral fixed point.  This map has no metric attractors.

\begin{coro}
The metric and topological attractors of $f$ coincide.
\end{coro}

{\it Decay of geometry.}
Following the concept of an adapted interval~\cite{mane} we  
call an open set $U$  \Def{regularly returning}
for some dynamics $f$ defined in an ambient space containing
$\overline{U}$ if $f^n(\partial U) \cap U = \emptyset$ for every  \pagebreak
$n > 0$. 

The \Def{first entry map} $E$ of $f$ into a set $U$ is defined
on $$E_U := \{x:\exists\,n > 0, f^{n}(x)\in U\}$$
by the formula $E(x) :=f^{n(x)}(x)$
where $n(x):=\min\{n>0: f^{n}(x)\in U\}$.
The \Def{first return map} of $f$ into $U$ is the restriction of the first
entry map to $E_U \cap U$.  The \Def{central domain} of the first return
map is the connected component of its domain containing the critical point
of $f$.
If $U$ is a regularly returning open interval then the function $n(x)$
is continuous and locally constant on $E_{U}$.

\begin{defi}\label{defi:24vp,1}
Suppose that  $J$ is an  open interval and $\overline{J} \subset I$. 
Define $\nu(J, I) := \frac{|J|}{\dist(J,\partial I)}$.
\end{defi}

The key property that enables us to exclude wild Cantor
attractors is the following result, known as \Def{decay of geometry}.

\begin{theo}\label{tw3}
Let $I$ be a compact interval and $f : I \to I$
be a $C^{3}$ unimodal map with $C^{3}$ nonflat critical point $\zeta$
of order $2$.
If $\zeta$ is recurrent and nonperiodic and
$f$ has only finitely many restrictive intervals then for every
$\eps_0 > 0$ there is a regularly returning interval
$Y' \ni \zeta$ such that if $Y$ is the central domain of the first
return map to $Y'${\rm ,} then $\nu(Y, Y') < \eps_0$.
\end{theo}

Decay of geometry occurs when the order of the
critical point is $2$.
Counterexamples exist when
the order of the critical point is larger than $2$~\cite{bkns}.

\demo{A priori bounds}
The following important fact known as
\Def{a priori bounds} is proved in~\cite[Lemma 7.4]{kozl}.
An earlier version for nonrenormalizable
maps can be found in~\cite{vargas}. 

\begin{fact}\label{fa:1}
Let $f$ be a $C^3$ unimodal map with $C^3$ nonflat
nonperiodic critical point $\zeta$.
Then there exists a constant $K$ and an
infinite sequence of pairs $Y'_{i}\supset Y_{i}\ni \zeta$
of open intervals with $|Y_i| \rightarrow 0$
such that{\rm ,} for each $i${\rm ,} $Y_i$ is regularly returning{\rm ,}
$\nu(Y_i, Y'_i) \leq K$ and for every branch $f^{n}$ of
the first entry map of $f$ into $Y_{i}${\rm ,} $f^n$
extends diffeomorphically
onto $Y'_{i}$ provided the domain of the branch is disjoint from $Y_{i}$.
\end{fact}

{\it Negative Schwarzian derivative and conjugation theorem.}
We say that a $C^3$ function $g$ has \Def{negative Schwarzian derivative}
if $$S(g)(x) := g'''(x)/g'(x) - \frac{3}{2} \left(g''(x)/g'(x)\right)^2 < 0$$
whenever $g'(x) \neq 0 $.  The Schwarzian derivative $S(g)$ satisfies
the composition law $S(g \circ h)(x) = S(g)(h(x)) h'(x)^2 +S(h)(x)$.
Thus iterates of a map with negative Schwarzian derivative also have
negative Schwarzian derivative.

In the general smooth case, negative Schwarzian derivative can be
recovered~\cite{comp} in the  \pagebreak following sense.

\begin{theo} \label{tw4}
Let $I$ be a compact interval and $f : I \to I$
be a $C^{3}$ unimodal map with $C^{3}$ nonflat and nonperiodic
critical point.  Then there exists a real\/{\rm -}\/analytic
diffeomorphism $h : I \to I$ and an {\rm (}\/arbitrarily small\/{\rm )}
open interval $U$ such that{\rm ,}
putting $g = h \circ f \circ h^{-1}${\rm ,} $U$ is a regularly returning
{\rm (}\/for $g${\rm )} neighborhood of the critical point of $g$ and the
first return map of $g$ to $U$ has uniformly negative Schwarzian
derivative.
\end{theo}  

 1.2. {\it  Box mappings}.
 \vglue-20pt
\phantom{badluck}

\begin{defi}\label{defi:704a,1}
Consider a finite sequence of compactly nested open
intervals around a point $\zeta\in
{\bf R}:\zeta\in b_0 \subset \overline{b}_0\subset b_1 \cdots \subset
b_k$. Let $\phi$ be a real-valued $C^1$ map defined on some open and bounded
set $U\subset {\bf R}$ containing $\zeta$.
Suppose that the derivative of $\phi$ only
vanishes at $\zeta$, which is a local extremum. Assume in addition the
following: 
\vglue-18pt
\phantom{way}
\begin{itemize} \itemsep=2pt
\item
for every $i=0,\cdots,k$, we have $\partial b_k \cap U = \emptyset$,
\item
$b_0$ is equal to the connected component of $U$ which contains $\zeta$,
\item
for every connected component $W$ of $U$ there exists $0 \leq i \leq k$
so that $\phi$ maps $W$ into $b_i$ and $\phi :\: W\rightarrow
b_i$ is proper.
\end{itemize}
\vglue-4pt\noindent
Then the map $\phi$ is called a {\em box mapping} and the intervals
$b_i$ are called {\em boxes}.
\end{defi}
\phantom{now}
\vglue-16pt

The restriction of a box map to a connected component of
its domain will be called a {\em branch}. Depending on
whether the domain of this branch contains the critical point $\zeta$ or not,
the branch will be called {\em folding} or {\em monotone}.
The domain $b_0$ of the folding branch is called the \Def{central domain}
and will usually be denoted by $b$;
$b'$ will denote the box into which the folding branch maps
properly. A box map $\phi$ is said to be {\em induced} by a map
$f$ if each branch of $\phi$ coincides on its domain
with an iterate of $f$ (the iterate may depend on the branch).

\vglue6pt {\it Type {\rm I} and type {\rm II} box mappings}.
A box mapping is of type II provided that there are only two boxes $b_0 = b$
and $b_1 = b'$, and every branch is proper into $b'$. A box mapping is of type
I if there are only two boxes, the folding branch is proper into $b'$
and all other branches are diffeomorphisms onto $b$.
A type I box mapping can be canonically obtained from a type II box map $\phi$
by {\em filling-in}, in which $\phi$ outside of $b$ is replaced by the first
entry map into $b$.
Note that if $f$ is a unimodal map with critical point $\zeta$ and $I$ is a
regularly returning open interval containing $\zeta$, then the first return
map of $f$ into $I$ is a type II box mapping.
 
\section{Distortion estimates}

In this section we prove a strong form of the $C^2$ Koebe lemma
(Proposition~\ref{LKoebe}).  In Lemma~\ref{lem:duncan} we
give a new proof of the required cross-ratio estimates.

Let $I$ be an open interval and
$h : I \to h(I) \subseteq \R$ be a $C^{1}$ diffeomorphism.
Let $a,b,c,d$ be distinct points of $I$ and
define the cross-ratio $\chi(a,b,c,d) :=\frac{(c-b) (d-a)}{(c-a)(d-b)}$.
By the distortion of $\chi$ by $h$ we mean
$$\kappa_{h}(a,b,c,d):=\chi(h(a),h(b),h(c),h(d))/\chi(a,b,c,d)\;.$$
We have the composition rule
\begin{equation}\label{comp}
\log \kappa_{g\circ h}(a,b,c,d)= \log \kappa_{g}(h(a),h(b),h(c),h(d))
+\log \kappa_{h}(a,b,c,d)\;.
\end{equation}

Define, for $x \neq y$,
$$K_h(x,y):=\frac{\partial}{\partial x} \log\left| \frac{h(x)-h(y)}{x-y}\right|
= \frac{h'(x)}{h(x)-h(y)} - \frac{1}{x-y} \;.$$
An elementary calculation shows that
$$\log \kappa_h(a,b,c,d) = \int_a^b K_h(x,c) - K_h(x,d) dx
=\int_{\partial R}K_{h}(x,y)dx\;.$$
where $R$ is the rectangle $[a,b] \times [c,d]$  
suitably oriented.
Note that $K_h(x,y)$ is perhaps integrated across the diagonal $x=y$.

We will also use $\rho_h(a,b,c,d) := \log \kappa_h(a,b,c,d) / (b-a)(d-c)$.

\begin{lem} \label{DLB}
Let $I$ be an open interval and let
$h : I \to h(I) \subseteq \R$ be a $C^2$ diffeomorphism
such that $1/\sqrt{|Dh|}$ is convex. Then
$\rho_h(a,b,c,d) \geq 0$ for all distinct points
$a,b,c,d$ in $I$.
\end{lem}

\Proof
If $h$ is $C^3$ then the result follows from the formula
$\log \kappa_h(a,b,c,d) = \int_a^b \int_c^d
1/(x-y)^2 - h'(x)h'(y)/(h(x)-h(y))^2dxdy$ since the integrand
is nonnegative (equivalent to a standard
inequality for maps with nonpositive Schwarzian derivative).
The $C^2$ statement follows by an approximation argument.
\hfill\qed

\begin{defi}
A continuous increasing function $\sigma : \R \to \R$
such that $\sigma(0)=0$ will be called a \Def{gauge function}.
\end{defi}

We first consider the case without critical points:
\begin{lem} \label{DLA}
Let $I$ be a compact interval and let $h : I \to h(I) \subseteq \R$
be a $C^2$ diffeomorphism.
Then there exists a gauge function $\sigma${\rm ,} for all distinct points $a, b, c, d$ in $I${\rm ,} such that
$\left|\rho_h(a,b,c,d)\right| \leq |\sigma(d-c)/(d-c)|$.
\end{lem}

{\it Proof}.
Extend $K_h(x,y)$ to the diagonal of $I \times I$ by
defining $K_h(x,x) = \frac{h''(x)}{2h'(x)}$ for $x \in I$. It is
easily checked using Taylor expansions that $K_h$ is continuous
and thus uniformly continuous.
Set $\Delta_{h}(x,y,z):=K_{h}(x,y)-K_{h}(x,z)$ and note that $\Delta_h(x,y,y) = 0$ for
all $x,y \in I$. Thus there exists
a gauge function $\sigma$ such that
$\left|\Delta_h(x,y,z)\right| \leq |\sigma(z-y)|$ for all $x,y,z \in I$.
From $\log \kappa_h(a,b,c,d) = \int_a^b \Delta_h(x,c,d)dx$ we see that
$\left|\log\kappa_h(a,b,c,d)\right| \leq |b-a|\,|\sigma(d-c)|$.
\hfill\qed
 
We now allow critical points on the boundary of the interval.
The following result generalizes a number of known cross-ratio
inequalities; see Theorems 2.1 and 2.2 of~\cite{S}.

\begin{lem}\label{lem:duncan}
Let $I$ be a compact interval and $f : I \to \R$ be a $C^2$ map with
all critical points $C^2$ nonflat.
Then there exists a gauge function $\sigma$
such that for all  distinct points $a,b,c,d$ in $I$ contained 
in the closure of a subinterval $J$ on which
$f$ is a diffeomorphism{\rm ,} 
$$ \rho_f(a,b,c,d) \geq -\min\left(\frac{\sigma(b-a)}{b-a}\;,\;
\frac{\sigma(d-c)}{d-c}\right)\;.$$
\end{lem}

\Proof 
It suffices to prove
$\rho_f(a,b,c,d) 
\geq - \sigma(d-c)/d-c$  since $\rho_f(a,b,c,d)=\rho_f(c,d,a,b)$. 
By $C^{2}$ nonflatness of the critical points,
for every  $c\in \Crit$   there exist
$\eps_{c}$ and a  $C^{2}$
diffeomorphism $\phi_{c}$ such that 
$f(x)=f(c)\pm |\phi_{c}(x)|^{\ell_{c}}$, $\ell_{c}>1$, 
 in $U_{c}=[c-\eps_{c},c+\eps_{c}]\cap I$.
Let $\eps:=\inf_{c\in \Crit} \eps_{c}/2$. 
Since $f$ has at most finitely many critical points, $\eps$
is positive. 

Suppose that $[a,d]$  is contained in an interval $J$ whose
endpoints are either in  $\Crit$ or in $\partial I$ and
$f'\neq 0$ inside $J$.  
Set $\Omega_{\eta} = \{(x,y) \in J^2 : \left|x-y\right| < \eta\}$
and note that
$K_f(x,y)$ is continuous for $(x,y)$ in the compact set
$J^{2} \setminus \Omega_{\eta}$.  
If  $[a,b] \times [c,d] \cap \Omega_{\eta} = \emptyset$ then 
\begin{equation}\label{far}
|\log \kappa_f(a,b,c,d)|= 
\left|\int_a^b K_f(x,c)-K_f(x,d) dx\right|
\leq |(b-a) \tilde{\sigma}(d-c)|
\end{equation}
for some gauge function $\tilde{\sigma}$ and consequently,
$|\rho_f(a,b,c,d)|\leq \tilde{\sigma}(d-c)/d-c$.

Now subdivide the rectangle $R=[a,b]\times [c,d]$
into $N$ equal rectangles $R_{i}=[a_{i},b_{i}]\times [c_{i},d_{i}]$ with 
the sides smaller than $\eta:=\eps/3$
and the orientation induced by $R$. 
In particular, the sign of
$(b_{i}-a_{i})(d_{i}-c_{i})$ does not depend on~$i$. We will use the fact that
\[\rho_{f}(a,b,c,d)=
\frac{1}{(b-a)(d-c)}\sum_{i}\int_{\partial R_{i}}K_{f}(x,y) dx\;=\;
 \frac{1}{N}\sum_{i}\rho_{f}(a_{i},b_{i},c_{i},d_{i})\;.\]
If $R_{i} \cap \Omega_{\eps/3}
 = \emptyset$ then the estimate $(\ref{far})$ works.  
If $R_{i} \cap \Omega_{\eps/3} \neq \emptyset$   
then $R_{i}$ is contained in 
$\Delta_{\eps}$. In particular, $a_{i},b_{i},c_{i},d_{i}$ are contained
in the interval $J_{i}$ of  length $\leq \eps$. We consider two cases.
 
\begin{description}   
\itemm{\rm (i)} If $J_{i}$ 
is not contained in $\bigcup_{c\in \Crit}U_{c}(\eps_{c})$
then the distance of $J_{i}$ to $\Crit$ is bigger than $\eps$.
To estimate $\int_{\partial R_{i}}K_{f}(x,y)dx$
 we apply Lemma~\ref{DLA} for $f$ restricted 
 $J\setminus \bigcup_{c\in \Crit}U_{c}(\eps)$.

\itemm{\rm (ii)}
If $J_{i}$ is  contained in $\bigcup_{c\in \Crit}U_{c}(\eps)$ then 
we write $f(x) = f(c) \pm |\phi_{c}(x)|^{\ell_{c}}$
for $x\in U_{c}(\eps_{c})$.
If $g=|\cdot|^\ell$ then $\rho_g(a_{i},b_{i},c_{i},d_{i})\geq 0$ and 
it is enough, by the composition rule
(\ref{comp}), to consider the effect of $\phi$. Lemma~\ref{DLA}
gives us the desired estimate.
\end{description}
We finish the proof 
by summing up the contributions  from all  rectangles $R_{i}$. \phantom{moresnow}
\Endproof

\begin{prop}[{\rm the Koebe principle}]\label{LKoebe}
Let $I$ be a compact interval and $f : I \to I$ be a $C^2$ map with
all critical points $C^{2}$ nonflat.
Then there exists a gauge function
$\sigma$ with the following property.
If $J\subset T$ are open intervals and $n\in N$ is such that
$f^{n}$ is a diffeomorphism on $T$ then{\rm ,} for every $x,y \in J${\rm ,} we have
$$
\frac{(f^n)'(x)}{(f^n)'(y)} \geq
\frac{e^{-\sigma(\max_{i=0}^{n-1} |f^i(T)|)
\sum_{i=0}^{n-1} |f^i(J)|}}
{(1 + \nu(f^n(J),f^n(T)))^{2}}\;.$$
\end{prop}

\Proof
Without loss of generality $T = (\alpha,\beta)$, $J = (x,y)$
and $\alpha < x < y < \beta$. Write $F = f^n$ and let $\sigma$ be as in
Proposition~\ref{lem:duncan}. Set
$P = \sum_{i=0}^{n-1} \sigma(|f^i(T)|) |f^i(J)|$.
By Proposition~\ref{lem:duncan} and (\ref{comp}),
\begin{eqnarray*}
\log\kappa_{F}(\alpha,x,x+\eps,y) &\geq&
\sum_{i=0}^{n-1}\log\kappa_{f}(f^{i}(\alpha),f^{i}(x),
f^{i}(x+\eps),f^{i}(y))\\
&\geq& - \sum_{i=0}^{n-1} \sigma(|f^i(\alpha,x)|) |f^i(x+\eps,y)|\\
&\geq& - \sum_{i=0}^{n-1} \sigma(|f^i(T)|) |f^i(J)| = -P.
\end{eqnarray*}
Taking $\eps \downarrow 0$ yields
$$\frac{F(y)-F(\alpha)}{y-\alpha}F'(x)
\geq e^{-P} \frac{F(x)-F(\alpha)}{x-\alpha}\frac{|F(J)|}{|J|}$$
which after rearranging becomes
\begin{equation} \label{KLE1}
F'(x)\geq e^{-P} \frac{|\alpha-y|}{|\alpha-x|}
\frac{|F(\alpha)-F(x)|}{|F(\alpha)-F(y)|}\;\frac{|F(J)|}{|J|}\geq
\frac{e^{-P}}{1+\nu(F(J),F(T))}\;\frac{|F(J)|}{|J|}\;.
\end{equation}
Likewise, considering $\kappa_F(x,y-\eps,x+\eps,y)$ and taking
$\eps \downarrow 0$ yields
$$\frac{|F(J)|^2}{|J|^2} \geq e^{-P} F'(x)F'(y)\;.$$
Equation (\ref{KLE1}) now
gives $F'(x)/F'(y)\geq e^{-3P}/(1+\nu(F(J),F(T)))^{2}\;$.
\hfill\qed

\section{Proof of the conjugation theorem}
In the following easy lemma we consider diffeomorphisms with constant
negative Schwarzian derivative.  These will be useful in defining the
conjugacy in Theorem~\ref{tw4}.

\begin{lem}\label{lem:25va,1}
For $s>0$ consider the function
\[
\phi_s(x) := \frac{\tanh(\sqrt{\frac{s}{2}}x)}{\tanh(\sqrt{\frac{s}{2}})} \; ,\]
which is a real\/{\rm -}\/analytic diffeomorphism of the real line into itself{\rm ,} fixing
$-1${\rm ,} $0$ and $1$.  The Schwarzian derivative of $\phi_s$ is everywhere equal
to $-s$.
\end{lem}

The following lemma is included for completeness.
An interval $J$ is \Def{symmetric} for a unimodal map $f$ if
$J = f^{-1}(f(J))$.

\begin{lem}\label{lem:sup}
Let $I$ be a compact interval and $f : I \to I$ be a unimodal map.
If $f$ does not have arbitrarily small regularly
returning symmetric open intervals containing the critical point $\zeta$ then
$\zeta$ is periodic.
\end{lem}

\Proof
Let $J$ be the interior of the intersection of all regularly returning
symmetric open intervals containing $\zeta$.  We must show that if
$J \neq \emptyset$ then $\zeta$ is periodic.
Indeed, if $J \neq \emptyset$ then $J$ is clearly a regularly returning
symmetric open interval containing $\zeta$.
By the minimality of $J$, $\zeta$ is mapped inside $J$ by some iterate
of $f$. Let $\phi$ be the first return map to $J$, which by minimality
has only one branch.
Again by minimality $\phi$ cannot have fixed points inside $J$ other than
$\zeta$.  Moreover $\zeta$ is indeed  a fixed point of $\phi$
since otherwise we could easily construct an appropriate regularly
returning interval inside $J$ containing~$\zeta$.
\Endproof

The next lemma is a standard consequence of the nonexistence
of wandering intervals \cite{dms}.

\begin{lem}\label{ops}
Let $f$ be a $C^{2}$ unimodal map with $C^2$ nonflat{\rm ,} nonperiodic
critical point $\zeta$.  For every interval $Y \ni \zeta$ there
exists $\eps_0(Y)>0$ such that if $I$ is an interval
mapped diffeomorphically onto $Y$ by some iterate $f^{n}$
then $|f^{i}(I)|\leq \eps_0(Y)$ for every $i=0,\cdots,n ${\rm ,} and
$\eps_0(Y) \to 0$ as $|Y| \to 0$.
\end{lem}

\Proof
Otherwise there exists $\delta>0$, a sequence
$Y_i \downarrow \{\zeta\}$ of open intervals, intervals $I_i$ with
$|I_i| > \delta$ and $n_i \to \infty$ such that $f^{n_i}$ maps
$I_i$ diffeomorphically onto $Y_i$.  Passing to a subsequence,
we may suppose the $I_i$ converge to some limit interval $I_\infty$ with
$|I_\infty| \geq \delta$.  Let $I$ be an interval compactly contained in the
interior of $I_\infty$.  By definition
$\Res{f^{n_i}}{I}$ is diffeomorphic for arbitrarily large $n_i$.  Thus 
$\Res{f^n}{I}$ is diffeomorphic for all $n > 0$, which shows that
$I$ is a homterval.  Since $f$ has no wandering intervals~\cite{dms},
this means that $\omega(x)$ is a periodic orbit for some $x \in I$.
However $\zeta \in \omega(x)$ by definition; thus $\zeta$ is periodic, a
contradiction.
\Endproof

Suppose now that $f^{n}$ is a branch of the first entry map into
an interval $Y:=Y_{i}$ given by fact~\ref{fa:1},
and that the domain $J$ of the branch is disjoint from $Y$.
There is a number $\eps(Y) > 0$ independent of the
branch such that for all $x \in J$ we have
$
S(f^n)(x) \leq \frac{\eps(Y)}{|J|^2}$ and
$\eps(Y) \to 0$ as $|Y| \to 0$.
Indeed, letting $L = \max(0,\sup S(f)) < \infty$, the composition law for
the Schwarzian derivative, the Koebe principle and the disjointness of
$J$, $\ldots$, $f^{n-1}(J)$ yield
\begin{eqnarray*}
S(f^n)(x) &=& \sum_{i=0}^{n-1} S(f)(f^i(x)) (f^i)'(x)^2
\leq \frac{L K^4}{|J|^2}\; \sum_{i=0}^{n-1} |f^i(J)|^2 \\
&\leq &
\frac{\eps_0 (Y) L K^4}{|J|^2} \; \sum_{i=0}^{n-1} |f^i(J)|
\leq
\frac{\eps (Y)}{|J|^{2}},
\end{eqnarray*}
where $K >0$ comes from the Koebe lemma.

We will also need the fact that if $f$ has a $C^3$
nonflat critical point $\zeta$
then there is some $\eta > 0$ such that
$S(f)(x) < - \eta / (x-\zeta)^2$ for all
$x \neq\zeta$ sufficiently close to $\zeta$.

\demo{Proof of Theorem~{\rm \ref{tw4}}}
Fix some $0 < s < \eta / 4$ and consider an interval
$Y:=Y_i$ given by fact~\ref{fa:1}.
Let $F$ denote the first return map to $Y$ and let $A$
be the increasing affine map taking $Y$ to $(-1,1)$.
Observe that $h := \phi_s \circ A$ has
constant Schwarzian derivative $S(h)(x) = -4s/ |Y|^2$.
The composition law for the Schwarzian derivative gives \setcounter{equation}{3}
\begin{equation}\label{equ: 1}
S(h\circ F\circ h^{-1})(h(x)) h'(x)^2 =
\frac{4s}{|Y|^2}(1-F'(x)^2) + SF(x)
\end{equation}
for all $x$ in the domain of $F$.
Let $I$ be the domain of a branch $f^{n+1}$ of $F$. Then $f(I)$ is contained
in a domain $J$ of a branch $f^n$ of the first entry map into $Y$, and $J$
is disjoint from $Y$ (if $\zeta \not\in I$ then
$J = f(I)$).
Let $G = h \circ F \circ h^{-1}$.
Equation $\ref{equ: 1}$ and the results noted above yield, for $x \in I$,
\begin{eqnarray*}
S(G)(h(x)) h'(x)^2 &=&
\frac{4s}{|Y|^2}(1- f'(x)^2 (f^n)'(f(x))^2) +
S(f^n)\circ f \; f'(x)^2 + S(f)(x)\\
&\leq & \frac{4s}{|Y|^2}\left(1- f'(x)^{2}\frac{|Y|^2}{K^2 |J|^2}\right)
+ \eps(Y)\frac{f'(x)^2}{|J|^2} + S(f)(x)\\
&=& \frac{f'(x)^2}{|J|^2} \left(\eps(Y) - \frac{4s}{K^2}\right)
+ S(f)(x) + \frac{4s}{|Y|^2}.
\end{eqnarray*}
Now $\eps(Y) - 4s/K^2$ will be negative as long as $Y$ is small enough.
Since $x \in Y$ we have $S(f)(x) + 4s/|Y|^2 < (4s-\eta)/|Y|^2 < 0$ if
$|Y|$ is small enough.
Thus $S(G)(y) < 0$ for all $y$ in the domain of $G$ if $Y$ is small enough.
\Endproof

We immediately obtain a weak form of the finiteness of attractors
theorem~\cite{dms}:

\begin{coro} \label{fott}
Let $I$ be a compact interval and $f : I \to I$ be a $C^3$ unimodal map
with $C^3$ nonflat critical point.  Then there exists $N \in \Z^+$ such
that any periodic orbit with period greater than $N$ is repelling.
\end{coro}

{\it Proof}.
It is well known~\cite[Th.\ C]{mane} that if
nonrepelling periodic orbits of arbitrarily high period exist, then
they must accumulate the critical point.  This is impossible if the critical
point $\zeta$ is periodic.  If $\zeta$ is not periodic then, by
the conjugation theorem, after a real-analytic coordinate change, the
first return map $\phi$ of $f$ to a regularly returning interval $U$
containing $\zeta$ has negative Schwarzian derivative.
Because of the negative Schwarzian derivative, all nonrepelling periodic
orbits of $\phi$ must attract $\zeta$,
so there can be at most one of these.  Since any periodic orbit of $f$ passing
through $U$ is periodic for $\phi$, this proves the result.
\hfill\qed

\section{Decay of geometry}

\begin{prop}\label{prop:26va,1}
Let $F$ denote a $C^{3}$ type {\rm II} box mapping with $C^{3}$ nonflat
critical point $\zeta$ of order $2$
and negative Schwarzian derivative. Assume that
the orbit of $\zeta$ is infinite{\rm ,} recurrent and that $F$ has no restrictive interval.
Suppose that $F$ has the following expansivity property\/{\rm :}
for every $\eta>0$ there is some $\delta>0$ such that
 if $I$ is an interval mapped by a nonnegative iterate
of $F$ diffeomorphically onto an interval of length less than $\delta$
containing $\zeta${\rm ,} then $|I| < ~\eta$.
Now{\rm ,} for every $\eps_0 > 0$ there is a regularly returning interval $Y'$
which contains $\zeta$ such that{\rm ,}
if $Y$ denotes the central domain of the first return map into $Y'${\rm ,}
then $\nu(Y,Y') < \eps_0$.
\end{prop}

The proof will be split into two cases depending on whether $\omega(\zeta)$ intersects the domains of finitely many
branches of $F$, or infinitely many.

\demo{The case with infinitely many branches}
For each domain $\Delta$ of a branch which intersects $\omega(\zeta)$,
we define $n(\Delta)$, the first entry time of the orbit of $\zeta$
into $\Delta$. We choose a sequence $\Delta_i$ with
$|\Delta_i| \rightarrow 0$. For each $i$ the interval $\Delta_i$ can be pulled
back by $F^{1-n(\Delta_i)}$ as a diffeomorphism to a neighborhood of the
critical value. The lengths of these neighborhoods tend to $0$ by the
expansivity hypothesis. Pulling back by one more iterate of $F$ we
get a sequence of regularly returning neighborhoods $Y'_i$ of $\zeta$, each of which is mapped into $\Delta_i$
by $F^{n(\Delta_i)}$ as a proper unimodal map.

Let us construct $Y_i$, the corresponding central domain.
It is the preimage by $F^{-n(\Delta_i)}$
of some domain $\delta_i$ of the first entry map into $Y'_i$.
We begin by considering the first entry map into the central
domain of $F$. Let $\delta'_i$ be the domain of the first entry map
into the central domain of $F$ which
contains $\delta_i$.
Since the nesting of the central domain inside the range of $F$ is preserved under pull-back by negative
Schwarzian diffeomorphisms,
we get $\nu(\delta'_i, \Delta_i) \geq \eta > 0$ where $\eta$ is independent of $i$.
Likewise, by the classical Koebe lemma for maps with negative Schwarzian derivative,
the distortion of the first entry map on $\delta'_i$ is bounded independently
of $i$. By the expansivity
hypothesis again, lengths of domains of branches of the first entry map into $Y'_i$ tend to $0$ with $i$.
Hence
$\nu(\delta_i, \Delta_i) \rightarrow 0$. Again by negative Schwarzian
derivative, this implies $\nu(Y_i, Y'_i) \rightarrow 0$ which implies the
result.

Observe that this argument did not make use of the hypothesis about the nondegeneracy of the critical point.
\vfill
\demo{The case with finitely many branches}
In this case will use the following theorem.
\vfill
\begin{theo}\label{theo:26vp,1}
Let $F$ denote a $C^3$ type {\rm II} box mapping with finitely many branches,
negative Schwarzian derivative{\rm ,} with the central branch which factors
as $F(\zeta)(1-(h(x))^2)$ where $h$ is a $C^3$
increasing diffeomorphism of the closure of the central domain onto 
its image and $h(\zeta) = 0$. Assume that the orbit of $\zeta$
is infinite{\rm ,} recurrent and that $F$ has no restrictive interval.
Now{\rm ,} for every $\eps_0 > 0$ there is a regularly returning interval
$Y'$ which contains $\zeta$ such that{\rm ,} if $Y$ denotes the central domain
of the first return map into $Y'${\rm ,} then $\nu(Y, Y') < \eps_0$.
\end{theo}
\vfill

Theorem~\ref{theo:26vp,1} is a weaker restatement of 
Theorem~1 from~\cite{nesz}. 
 
The map $F$ from Proposition~\ref{prop:26va,1} can be restricted to only
those branches whose domains intersect $\omega(\zeta)$ and the remaining
hypotheses will still apply.
Now, we see that Theorem~\ref{theo:26vp,1} is directly applicable under
the hypotheses of Proposition~\ref{prop:26va,1} except for the special form
of the central branch. This problem is taken care  of by an elementary
calculation.

\vfill{\it Proof of Theorem~{\rm \ref{tw3}}}.
We apply Proposition~\ref{prop:26va,1} to the first return map $F$ of $f$ to
a regularly returning open interval $U$ containing the critical point.
We may assume that $F$ has negative Schwarzian derivative by
Theorem~\ref{tw4}.  If $U$ is small enough then $F$ will have no
restrictive intervals.  Moreover, we take $U$ small enough that its
closure contains no nonrepelling periodic points (Corollary~\ref{fott}).
The expansivity hypothesis is then satisfied.  Indeed, if not, then
arguing by contradiction as in Lemma~\ref{ops}, we see that $F$ (and
thus $f$) would have a homterval $I$.  Since $f$ has no wandering intervals,
this means that some point of $I$ is attracted to a nonrepelling periodic
orbit.  This orbit must intersect the closure of $U$, a contradiction.
\Endproof
\vfill

{\it Induced expansion.}
We will say that a unimodal map $f$ \emph{induces expansion} if there is
a regularly returning open interval $J$ containing the critical point of $f$,
an open subset $U$ of $J$, a map $F : U \to J$ and  $\rho > 1$ with the
following properties: for each connected component $V$ of $U$ there is a
positive integer $n(V)$ such that $F$ coincides with $f^{n(V)}$ on $V$,
$F$ maps $V$ diffeomorphically onto $J$ with derivative (in absolute value)
at least $\rho$, and if $A$ is the set of points in $J$ which return to
$J$ infinitely often under iteration by $f$, then $U$ contains Lebesgue almost
every point of $A$.\eject

\begin{prop} \label{cie}
Let $I$ be a compact interval and $f : I \to I$
be a $C^{3}$ unimodal map with $C^{3}$ nonflat critical point $\zeta$
of order $2$.
If $\zeta$ is nonperiodic and $f$ has only finitely many restrictive
intervals then $f$ induces expansion.
\end{prop}

{\it Proof}.
As in the proof of Theorem~\ref{tw3}, we consider the first return
map $F$ of $f$ to a regularly returning interval $U$ containing the
critical point.  We may assume that $F$ has negative Schwarzian
by Theorem~\ref{tw4}.  By Theorem~\ref{tw3}, we may take the central
domain of $F$ to be as small as we like, proportionally to $U$, by
taking $U$ small enough.  Adapting Proposition 5 from~\cite{jas} to
$F$, we see that $F$ induces an expanding Markov map on some perhaps
smaller regularly returning open interval $U'$ containing $\zeta$,
in other words $f$ induces expansion on $U'$.
\hfill\qed

\section{Attractors}

{\it Dynamics away from critical points.}
Our study of the dynamics away from critical points is based on the
following result of Ma{\~n}\'e~\cite[Th.\ D]{mane}:

\begin{fact}\label{pmane}
Let $I$ be a compact interval{\rm ,} $f : I \to I$ be a
$C^2$ map and $K \subseteq~I$ be a compact invariant set not containing
critical points.  Then either $K$ has Lebesgue measure zero or there
exist an interval $J \subseteq I$ and $n \geq 1$ such that
$f^n(J) \subseteq J${\rm ,} $\Res{f^n}{J}$ has no critical points and
$J \cap K$ has positive Lebesgue measure.
\end{fact}

An interval map $f$ is \Def{nonsingular} if $f^{-1}(A)$ has
zero Lebesgue measure for every Borel set $A$ with zero Lebesgue measure.
A map with a finite number of critical points is nonsingular.

\begin{coro}\label{nonhyp}
Let $I$ be a compact interval and $f : I \to I$ be $C^2$
and nonsingular.  Then{\rm ,} for Lebesgue almost every $x\in I${\rm ,}
either $\omega(x)$ contains a critical point or $\omega(x)$ coincides with
a nonrepelling periodic orbit.
\end{coro}

{\it Proof}.
Let $A$ be the set of all $x \in I$
such that  $\omega(x) \cap \Crit = \emptyset$ and $\omega(x)$
is not a nonrepelling periodic orbit.  If $A$ is of positive Lebesgue measure
then there exists an open neighborhood $U$ of $\Crit$ and a forward
invariant $B\subset A$ 
of positive Lebesgue measure so that the forward orbit
of every  point of $B$ is disjoint from $U$. The support $K$
of Lebesgue measure restricted to $B$ 
is forward invariant also.  By fact~\ref{pmane},
there exist  an interval $J \subseteq I$ and $n \geq 1$ such that
$f^n(J) \subseteq J$, $\Res{f^n}{J}$ has no critical points and
$J \cap K$ has positive measure.  It follows that $J \cap A$ has
positive measure also.  But $\omega(x)$
is a nonrepelling periodic orbit for almost every $x \in J$, a contradiction.
\hfill\qed \eject

\begin{coro}\label{bwump}
Let $I$ be a compact interval and $f : I \to I$ be $C^2$
and nonsingular.  Then every metric attractor of $f$ that does not contain
critical points coincides with a topologically attracting periodic orbit.
\end{coro}

{\it Proof}.
Let $K$ be a metric attractor that contains no critical points.  Since
$\omega(x) = K$ for almost every $x \in \rho(K)$, it follows from
the preceding corollary that $K$ is a nonrepelling periodic orbit.
If this orbit is not topologically attracting then $\rho(K)$
coincides with $\boldsymbol{\cup}_{n \geq 0} f^{-n}(K)$ and has
measure zero, a contradiction.
\hfill\qed
\vglue12pt

It may be instructive to examine Figure~1 in the light of these
results.

\demo{Attractors containing the critical point}
In the light of Corollary~\ref{nonhyp}, Theorem~\ref{tw2} is reduced to
the following assertion.

\begin{prop}\label{lem:a0}
Let $I$ be a compact interval and $f : I \to I$ be a $C^3$
unimodal map with $C^3$ nonflat critical point of order $2$.  Then{\rm ,}
for Lebesgue almost every $x \in I${\rm ,} either $\omega(x)$ does not contain
a critical point or $\omega(x)$ coincides with either
\begin{itemize}
\item a super\/{\rm -}\/attracting periodic orbit{\rm ,} or
\item a transitive cycle of intervals{\rm ,} or
\item a Cantor set of solenoid type.
\end{itemize}
\end{prop}

The proof of the proposition uses the following lemma.
\begin{lem}~\label{interval}
Let $I$ be a compact interval and $f : I \to I$ be a unimodal map
with critical point $\zeta$.  If there exists a subinterval $Y$ containing
$\zeta$ in its interior such that{\rm ,} for Lebesgue almost every $x \in Y${\rm ,}
the orbit of $x$ intersects $Y$ in a set of full Lebesgue measure{\rm ,} then
$f$ has a metric attractor which is a transitive cycle of intervals.
\end{lem}

{\it Proof}.
Note that some iterate of $Y$ intersects $Y$ since almost every point in $Y$
returns to $Y$.  Thus $C$, the closure of the union of the iterates of $Y$, is
a finite union of intervals.  Clearly $C$ is forward invariant and contains
$\zeta$ in its interior.  By the definition of $Y$, the orbit of almost every
$x \in C$ is dense in $C$ which implies that $C$ is a transitive metric
attractor.
\Endproof

{\it Proof of Proposition}~\ref{lem:a0}.
If the critical point is periodic then it is a super-attracting periodic
point and the result is obvious.  If $f$ is infinitely renormalizable
and $\zeta \in \omega(x)$, then $\omega(x)$ is obviously contained in the
intersection of the orbits of all restrictive intervals.  It is easily shown
that this is a Cantor set of solenoid type \pagebreak coinciding with $\omega(c)$.
In fact, one can even use the
classical\break S-unimodal theory since, by Theorem~\ref{tw4}, after a real-analytic
coordinate change, all renormalizations of $f$ of sufficiently high
period have negative Schwarzian derivative.

We may therefore suppose that the critical point is not recurrent and
that $f$ has only finitely many restrictive intervals.
From Corollary~\ref{cie} we know that $f$ induces expansion on an open
interval $J$ containing the critical point,
the induced map $F$ being defined on an open subset $U$ of $J$.
If $U$ has full (Lebesgue) measure in $J$ then the orbit by $F$ of
almost every point in $J$ is dense in $J$; thus, by  Lemma~\ref{interval}, $f$ has a metric attractor which
is a transitive cycle of intervals.  If $U$ does not
have full measure in $J$ then almost every
point in $U$ leaves $U$ under iteration by $F$.  By the definition
of induced expansion, almost every point $x \in J$ with $c \in \omega(x)$
never leaves $U$ under iteration by $F$.  Thus in this case we see that
the set of points $x$ with $c \in \omega(x)$ has measure zero, and  so the
result is trivial.
  \Endproof

{\it Proof of Corollary}~\ref{tw2}.
Suppose that $A$ is a metric attractor of $f$. If $A$
contains the critical point then, by Proposition~\ref{lem:a0}, $A$ is either
a super-attracting periodic orbit or a cycle of intervals or a solenoid.
It follows from their definitions that these three possibilities are
mutually exclusive.
If $A$ does not contain the critical point then,
by Corollary~\ref{bwump}, $A$ is a topologically
attracting periodic orbit.
\Endproof

\end{document}